\documentclass[a4paper,12pt]{amsart}

\usepackage{t1enc}
\def\R{\mathbb{R}}

\def\cE{\mathcal{E}}

\def\al{\alpha}

\def\de{\delta}

\def\si{\sigma}

\def\om{\omega}

\def\Om{\Omega}
\def\Up{\Upsilon}

\def\P{\mathrm{P}}

\newtheorem{theorem}{Theorem}[section]

\newtheorem{corollary}[theorem]{Corollary}
\theoremstyle{remark}
\newtheorem{remark}[theorem]{\rm\bf Remark}
\newtheorem*{remark*}{\rm\bf Remark}

\usepackage{amssymb,stmaryrd}
\usepackage{amscd}

\newcommand{\nd}{\nabla}

\makeatletter
\newcommand{\semidownbracefill}{$\m@th\braceld\leaders\vrule\hfill\braceru
  \bracelu\leaders\vrule\hfill\arrowhead$}



\def\sideremark#1{\ifvmode\leavevmode\fi\vadjust{\vbox to0pt{\vss
 \hbox to 0pt{\hskip\hsize\hskip1em
 \vbox{\hsize3cm\tiny\raggedright\pretolerance10000
 \noindent #1\hfill}\hss}\vbox to8pt{\vfil}\vss}}}%
                        


   \thanks{The author would like to thank Mike Eastwood and the referee for comments and suggestions.}

   \dedicatory{Dedicated to Mike Eastwood on the occasion of his 60th birthday}

\author{Matthew Randall}
 \address{Mathematical Sciences Institute, Australian National University, Canberra, ACT 0200, Australia}
\email{matthew.randall@anu.edu.au}
\title{Local obstructions to a conformally invariant equation on M\"obius surfaces}

\begin{document}

\begin{abstract}
On a M\"obius surface, as defined in \cite{mob}, we study a variant of the Einstein-Weyl (EW) equation which we call scalar-flat M\"obius EW (sf-MEW). This is a conformally invariant, finite type, overdetermined system of semi-linear partial differential equations. We derive local algebraic constraints for this equation to admit a solution and give local obstructions. In the generic case when a certain invariant of the M\"obius structure given by a symmetric tensor $M_{ab}$ is non-zero, the obstructions are given by resultants of $3$ polynomial equations whose coefficients are conformal invariants of the M\"obius structure. The vanishing of the resultants is a necessary condition for there to be solutions to sf-MEW.      
\end{abstract}

\maketitle

\pagestyle{myheadings}
\markboth{Randall}{Local obstructions to a conformally invariant equation on M\"obius surfaces}

\section{Introduction}
Let $(M^2,[g])$ be a Riemann surface. This is a smooth $2$-dimensional oriented manifold equipped with a conformal structure $[g]$, which is an equivalence class of smooth Riemannian metrics under the equivalence relation $g_{ab}\mapsto \widehat g_{ab}=\Om^2 g_{ab}$ for any smooth positive nowhere vanishing function $\Om$. Since every metric $g_{ab}$ in dimension $2$ is locally a conformal rescaling of the flat metric $\de_{ab}$, conformal geometry in dimension $2$ carries no local information. To remedy this, one can impose on Riemann surfaces additional local structure present in conformal manifolds of dimension $n>2$. This is the motivation behind M\"obius structures \cite{mob}. A Riemann surface with a M\"obius structure will henceforth be called a M\"obius surface. We can study on M\"obius surfaces a well-defined conformally invariant equation, which we call the scalar-flat M\"obius Einstein-Weyl (sf-MEW) equation. This equation specialises the Einstein-Weyl equation in conformal geometry in higher dimensions to the $2$-dimensional setting. We derive algebraic constraints for a given M\"obius surface to admit a solution to sf-MEW, and from there derive obstructions to existence of solutions to the equation. Checking that the obstructions do not vanish tells us definitively that the M\"obius surface  cannot admit any solution to the sf-MEW equation locally. Abstract indices \cite{pr} will be used throughout the paper to describe tensors on the conformal manifold. We have already used $g_{ab}$ to denote the metric tensor. For another instance, if we write $\om_{a}$ to denote a smooth $1$-form $\om$, then the $2$-form $d\om$ can be written as $\nd_{[a}\om_{b]}=\frac{1}{2}(\nd_a\om_b-\nd_b\om_a)$.   

\section{Conformal geometry and M\"obius structures}

A M\"obius surface is a Riemann surface $(M^2,[g])$ equipped with a smooth M\"obius structure as defined in $2.1$ of \cite{mob}. Taking the Weyl derivative $D_a$ in the definition to be the Levi-Civita connection $\nd_a$ for a particular representative metric $g_{ab}$ in the conformal class $[g]$, a M\"obius structure on $(M^2,[g])$ is a smooth second order linear differential operator $D_{(ab)_\circ}: \cE[1] \to \cE_{(ab)_\circ}[1]$ such that $D_{(ab)_\circ}-\nd_{(a}\nd_{b)_\circ}$ is a zero order operator acting on sections of the density line bundle of weight $1$, denoted by $\cE[1]$. Let $\P_{(ab)_\circ}$ be the symmetric trace-free tensor denoting the difference, i.e.\
\begin{align*}
\P_{(ab)_\circ}\si:=(D_{(ab)_\circ}-\nd_{(a}\nd_{b)_\circ})\si,
\end{align*}  
where $\si$ is a section of $\cE[1]$. Since the operator $D_{(ab)_\circ}$ is invariantly defined, under a conformal rescaling of the metric $\widehat g_{ab}=\Om^2g_{ab}$ we find that 
\begin{align*}
\widehat \P_{(ab)_\circ}=\P_{(ab)_\circ}-\nd_a\Up_b+\Up_a\Up_b-\frac{1}{2}g_{ab}\Up_c\Up^c+\frac{1}{2}g_{ab}\nd_c\Up^c 
\end{align*} 
where $\Up_a=\nd_a\log \Om$. Let $K$ denote the Gauss curvature of $g_{ab}$, i.e.\ $K=\frac{R}{2}$, where $R$ is the scalar curvature of $g_{ab}$. Define the Rho tensor by $\P_{ab}:=\P_{(ab)_\circ}+\frac{K}{2}g_{ab}$. Under a conformal rescaling, $K$ transforms as $\widehat K=K-\nd_a\Up^a$ and therefore 
\begin{align*}
\widehat \P_{ab}=\P_{ab}-\nd_a\Up_b+\Up_a\Up_b-\frac{1}{2}g_{ab}\Up_c\Up^c.
\end{align*} 
Hence for any representative metric $g_{ab}$ in the conformal class $[g]$ with its associated Levi-Civita connection $\nd_a$, a M\"obius structure in the sense of \cite{mob} determines a symmetric tensor $\P_{ab}$ satisfying the following two properties: \\
1) The metric trace of $\P_{ab}$ is the Gauss curvature $K$ of $g_{ab}$; \\
2) Under a conformal rescaling of the metric $\hat g_{ab}=\Om^2g_{ab}$, the tensor $\P_{ab}$ transforms accordingly as
\begin{equation}\label{Rho}
\widehat \P_{ab}=\P_{ab}-\nd_a\Up_b+\Up_a\Up_b-\frac{1}{2}g_{ab}\Up_c\Up^c,
\end{equation}     
where $\Up_a=\nd_a\log \Om$. This is the definition of a M\"obius structure given in \cite{mobreview} which we shall subsequently use. A M\"obius surface will be given by $(M^2,[g],[\P])$, where $[\P]$ denotes the M\"obius structure on $(M^2,[g])$, which is the equivalence class of smooth symmetric tensors related to a representative Rho tensor $\P_{ab}$ by formula (\ref{Rho}) under conformal rescalings of the metric. In $2$ dimensions, the Schouten tensor $\P_{ab}$ is not well-defined and a M\"obius structure remedies that by equipping the manifold with a Rho tensor that behaves like the Schouten tensor under conformal rescalings. In $2$ dimensions, a Rho tensor $\P_{ab}$ allows us to write 
\[
R_{abcd}=K(g_{ac}g_{bd}-g_{bc}g_{ad})\equiv \P_{ac}g_{bd}-\P_{bc}g_{ad}+\P_{bd}g_{ac}-\P_{ad}g_{bc},
\]
even though the tensor $\P_{ab}$ cannot be recovered from the Riemannian curvature tensor alone, in contrast to the higher dimensional setting. A fixed representative metric $g_{ab}$ from the conformal class $[g]$ can be viewed as having conformal weight $2$ and induces a volume form $\epsilon_{ab}=\epsilon_{[ab]}$ of conformal weight $2$. We set our convention so that $\epsilon^{ab}\epsilon_{cb}=\de_c{}^a$ and we raise and lower indices using the metric. The Cotton-York tensor given by
\[
Y_{abc}=\nd_a\P_{bc}-\nd_b\P_{ac}
\]
is a conformal invariant of the M\"obius structure. This means that under conformal rescalings of the metric, the quantity $Y_{abc}$ remains unchanged. 
We can use the volume form $\epsilon_{ab}$ to dualise, so that
\[
Y_{abc}=\frac{1}{2}\epsilon_{ab}Y_c, 
\]
where $Y_c=\epsilon^{ab}Y_{abc}$
is now a $1$-form of conformal weight $-2$. 
Observe that
\[
\frac{1}{2}\epsilon_{ab}Y^b=Y_{ab}{}^b=\nd_aK-\nd^b\P_{ab}.
\]
The vanishing of $Y_a$ characterises flat M\"obius surfaces. 

\section{The sf-MEW equation on M\"obius surfaces}
A Weyl connection $D_a$ on a Riemann surface $(M^2,[g])$ is a torsion-free connection that preserves the conformal class $[g]$, or equivalently $D_a g_{bc}=2\al_ag_{bc}$ for some $1$-form $\al_a$. The $1$-form $\al_a$ is determined up to a gauge freedom; under a conformal rescaling of the metric $g_{ab}\mapsto \widehat g_{ab}=\Om^2 g_{ab}$, we have $\al_a \mapsto \widehat \al_a=\al_a+\Up_a$, where again $\Up_a=\nd_a\log \Om$. For a fixed M\"obius surface $(M^2,[g],[\P])$, we can ask whether there is a compatible Weyl connection $D_a$ such that the second order linear differential operator $D_{(ab)_\circ}$ is given by the trace-free symmetric Hessian of the Weyl derivative, i.e.\ $D_{(ab)_\circ}=D_{(a}D_{b)_\circ}$. For a representative metric $g_{ab} \in [g]$ and its associated Levi-Civita connection $\nd_a$, this is equivalent to solving the system of equations given by
\begin{equation}\label{mew}
\mbox{Trace-free part of }(\nd_{(a}\al_{b)}+\al_a\al_b+\P_{ab})=0.
\end{equation}
Equation (\ref{mew}) is not finite type in the sense of \cite{spencer} because of its symbol. (We can compare this to the conformal Killing equation, which also has the same symbol and is also not finite type in 2 dimensions). However, if we impose the additional condition that the scalar curvature of the Weyl connection $D_a$ is $0$, the equation becomes finite type. This is the sf-MEW equation and is given by
\begin{equation}\label{sfmew0}
\nd_{(a}\al_{b)}+\al_a\al_b+\P_{ab}-\frac{\al_c\al^c}{2}g_{ab}=0.
\end{equation}
The `scalar-flat' in the terminology of sf-MEW hence refers to the vanishing of the scalar curvature of $D_a$. The sf-MEW equation is a conformally invariant, finite type, overdetermined system of semi-linear partial differential equations. A common procedure to treat equations such as (\ref{sfmew0}) is through prolongation \cite{prolong}. This involves expressing first derivatives of the dependent variables in terms of the variables themselves. Let $F_{ab}=\nd_{[a}\al_{b]}=\frac{1}{2}\epsilon_{ab}F$ be the extra dependent variable where $F=\epsilon^{ab}F_{ab}$. We can rewrite (\ref{sfmew0}) as
\begin{equation}\label{sfMEW}
\nd_a\al_b+\al_a\al_b+\P_{ab}-\frac{\al^c\al_c}{2}g_{ab}=\frac{1}{2}\epsilon_{ab}F.
\end{equation}
Observe that tracing the equation gives $\nd_a\al^a+K=0$ as a necessary constraint of the system.  
Differentiating (\ref{sfMEW}), we find that
\[
\nd_aF=-2\al_aF-Y_a
\]
is a consequence of the original equation. 
The derivative of the extra dependent variable $F$ is now given by known quantities $\al_a$, $F$ and $Y_a$ of the system. 
The prolonged system is therefore given by
\begin{align}
\nd_a\al_b=&\frac{1}{2}\epsilon_{ab}F+\frac{\al^c\al_c}{2}g_{ab}-\al_a\al_b-\P_{ab}\label{mewf1},\\
\nd_aF=&-2\al_aF-Y_a\label{mewf2}.
\end{align}
We can use the prolonged system to derive algebraic constraints for there to be a solution to (\ref{sfmew0}). We obtain as a necessary condition for there to be solutions of (\ref{sfmew0}) three polynomial equations in a single variable $t$ with coefficients given by conformal invariants of the M\"obius structure, under the assumption that a certain conformal invariant $M_{ab}$ defined in (\ref{defm}) is non-zero. This is the result of Theorem \ref{polycon}. The resultants of any two of these $3$ polynomials will then have to vanish for there to be a common root $t$, and since the resultants are given only by the invariants of the M\"obius structure, we obtain obstructions for there to be solutions of (\ref{sfmew0}) in Corollary \ref{polyobs}. We also discuss the case where $M_{ab}=0$ in Theorem \ref{mab}. We conclude the paper in section \ref{examples} by giving three examples of M\"obius structures on $\R^2$ with the flat metric $\de_{ab}$ for which one admits a solution to (\ref{sfmew0}), one does not because of non-vanishing obstructions and one with vanishing obstructions but does not admit any real solution to (\ref{sfmew0}). 
\begin{remark}
We remark on the case of global solutions to (\ref{sfmew0}) on closed M\"obius surfaces $(M^2,[g],[\P])$. Since $\nd_a\al^a+K=0$ is a consequence, we find that integrating this equation over $M^2$, the integral of the divergence term vanishes, so that $0=-\int_{M^2} K=-2\pi\chi(M^2)$ by the Gauss-Bonnet formula and this implies that $M^2$ has to be the torus.  
\end{remark}   

\section{Deriving algebraic constraints for (\ref{sfMEW}) to hold}\label{polyconproof}
In this section we derive the algebraic constraints that have to be satisfied for equation (\ref{sfMEW}) to hold. In the flat case when $Y_a=0$, we necessarily have $F=0$ by differentiating (\ref{mewf2}) and skewing. The $1$-form $\al_a$ is therefore exact, and equation (\ref{sfMEW}) specialises to the conformally Einstein equation for flat M\"obius surfaces. We therefore restrict our attention to non-flat M\"obius structures, that is one with $Y_a$ non-zero. By (\ref{mewf2}) this ensures that $F\neq 0$.   
Introduce the vector $U^a=\epsilon^{ab}Y_b$. Locally this is obtained by rotating the vector $Y^a$ $90^{\circ}$ clockwise on the plane. We have $Y_b=U^a\epsilon_{ab}$. The $1$-form $U_a$ has conformal weight $-2$. 
Differentiating (\ref{mewf2}) and skewing with $\epsilon^{ac}$ gives the first constraint of the system, namely that
\begin{equation}\label{mewcon1}
-2F^2=\epsilon^{ac}\nd_aY_c+2\epsilon^{ac}Y_c\al_a=\nd_aU^a+2\al_aU^a.
\end{equation}
Define the quantities
\begin{align*}
\mu:=\frac{\nd_cY^c}{2},\qquad \phi:=\frac{\nd_cU^c}{2}.
\end{align*} 
We can rewrite equation (\ref{mewcon1}) as
\begin{equation}\label{mewcon2}
\al_aU^a=-F^2-\phi.
\end{equation}
Set $W_a= Y^c\nd_cU_a+\phi Y_a-3 \mu U_a$.
The 1-form $W_a$ can be verified to be conformally invariant of weight $-6$. 
Differentiating (\ref{mewcon2}) and using (\ref{mewf1}), we find that 
\begin{align*}
\al_aW^a
=&3(\al_bY^b)F^2+\frac{5}{2}\rho F-Y^b\nd_b\phi+\P_{ba}U^aY^b+3\mu F^2+3\mu\phi,
\end{align*}
where $\rho:=U_aU^a=Y_aY^a$.
To simplify notation, set
\begin{align*}
\ell=&3\mu \phi+\P_{ab}U^aY^b-Y^c\nd_c\phi.
\end{align*}
We have
\begin{equation}\label{alphaw}
\al_aW^a=\ell+\frac{5}{2}\rho F+(3\mu+3\al_aY^a)F^2.
\end{equation}
We can now solve for $\al_a$ assuming $\si:=Y_aW^a \neq 0$ ($\si$ has conformal weight $-10$). It is given by
\begin{align}\label{alpha}
\al_a
=&\frac{Y_a}{\si}\left(\ell+\frac{5}{2}\rho F+(3\mu+3\al_cY^c)F^2\right)-\frac{\epsilon_{ab}W^b}{\si}\left(F^2+\phi\right).
\end{align} 
Contracting with $Y^a$ on both sides gives
\begin{align*}
\al_aY^a
=&\frac{\rho}{\si}\left(\ell+\frac{5}{2}\rho F+(3\mu+3\al_cY^c)F^2\right)+\frac{\tau}{\si}\left(F^2+\phi\right),
\end{align*}
where $\tau:=U_aW^a$ has conformal weight $-10$.  
Hence
\begin{equation}\label{abs}
\left(1-\frac{3\rho F^2}{\si}\right)\al_aY^a=\frac{\rho}{\si}\left(\ell+\frac{5}{2}\rho F+3\mu F^2\right)+\frac{\tau}{\si}\left(F^2+\phi\right), 
\end{equation}
and for $P_0(F):=\si-3\rho F^2\neq 0$, we obtain 
\begin{align*}
\al_aY^a=&\frac{\rho}{\si-3\rho F^2}\left(\ell+\frac{5}{2}\rho F+3\mu F^2\right)+\frac{\tau}{\si-3\rho F^2}\left(F^2+\phi\right).
\end{align*}
Substituting this expression into (\ref{alpha}) gives
\begin{align*}
\al_a=&\frac{Y_a}{\si}\left(\ell+\frac{5}{2}\rho F+3 \mu F^2\right)-\frac{\epsilon_{ab}W^b}{\si}\left(F^2+\phi\right)\\
&+\frac{3Y_a\rho F^2}{\si(\si-3\rho F^2)}\left(\ell+\frac{5}{2}\rho F+3 \mu F^2\right)+\frac{3Y_a\tau F^2}{\si(\si-3\rho F^2)}\left(F^2+\phi\right).
\end{align*}
which simplifies to 
\begin{align*}
\al_a=
&\frac{Y_a}{\si-3\rho F^2}\left(\ell+\frac{5}{2}\rho F+3\mu F^2\right)-\frac{\epsilon_{ab}W^b-3F^2U_a}{\si-3\rho F^2}\left(F^2+\phi\right).
\end{align*}
We can rewrite this expression to obtain
\begin{align*}
(\si-3\rho F^2)\al_a=&\left(Y_a\ell-\epsilon_{ab}W^b\phi\right)+\frac{5}{2}Y_a\rho F+\left(\frac{\nd_a\rho}{2}\right)F^2+3F^4U_a.
\end{align*}
The case where $P_0(F)=\si-3\rho F^2=0$ will be discussed in section \ref{mabproof}.
Call $L_a=Y_a\ell-\epsilon_{ab}W^b\phi$.
The 1-form $L_a$ has conformal weight $-10$, and we have
\begin{equation}\label{alpha1}
(\si-3\rho F^2)\al_a=L_a+\frac{5}{2}Y_a\rho F+\frac{\nd_a\rho}{2}F^2+3F^4U_a.
\end{equation}
The derivation of (\ref{alpha1}) still holds when $\si=0$. We shall now derive polynomial constraints for there to admit a solution of (\ref{sfMEW}). This involves differentiating equation (\ref{alpha1}) and using the equations (\ref{mewf1}) and (\ref{mewf2}) in the prolonged system to substitute. We then contract by the quantities $U^aU^b$, $Y^aY^b$ and $\epsilon^{ab}$ to produce three polynomial constraint equations 
\[
P_1(F)=0, P_2(F)=0, P_3(F)=0,
\]
in the variable $F$ with coefficients given by conformal invariants of the M\"obius structure that have to be satisfied for there to be solutions to (\ref{sfMEW}). After a routine computation we find that the first polynomial constraint is given by
\begin{align*}
P_1(F)=&\frac{63}{2}\rho^2 F^8-12\rho \si F^6+\begin{pmatrix}12\rho \si \phi-63\rho^2\phi^2+3\rho U^a\nd_a\si +\frac{1}{2}(\tau+3\mu \rho)^2\\ +\frac{1}{2}(3\rho \phi-\si)^2+\frac{3}{2}\rho U^aU^b\nd_a\nd_b\rho-9\rho^2\P_{ab}U^aU^b\end{pmatrix}F^4\\
&+\left(\frac{15}{2}\rho^3 \mu+\frac{5}{2}\tau \rho^2+\frac{15}{2}\rho^2(U^aU^b\nd_bY_a)\right)F^3\\
&+\begin{pmatrix}(3\rho\phi-\si)(U^a\nd_a\si)+21\rho \phi^2\si-3\phi \si^2+(\rho\ell+\phi\tau)(3\mu\rho+\tau)\\+\frac{25}{8}\rho^4+3\rho U^aU^b\nd_bL_a+6\rho\si \P_{ab}U^aU^b-\si\frac{U^aU^b\nd_a\nd_b\rho}{2}\end{pmatrix}F^2\\
&+\left(\frac{5}{2}\rho^2(\rho\ell+\phi \tau)-\frac{5}{2}(U^aU^b\nd_bY_a)\si \rho\right)F\\
&-\si\phi(U^a\nd_a\si)+\frac{1}{2}(\rho\ell+\phi \tau)^2-\frac{1}{2}\phi^2\si^2-\si(U^aU^b\nd_bL_a+\si \P_{ab}U^aU^b)\\
=&0.
\end{align*}
The second constraint is given by
\begin{align*}
\tilde P_2(F)=&\left(\frac{\si-15 \rho F^2}{\si-3\rho F^2}\right)\left((\rho\ell+\tau \phi)+\frac{5}{2}\rho^2 F+(\tau+3\mu \rho)F^2\right)^2-\frac{9}{2}\rho^2F^8\\
&-(9(Y^bY^a\nd_bU_a)\rho+3\rho(3\phi \rho-\si))F^6-25\rho^2(\tau+3\mu\rho)F^3\\
&+\begin{pmatrix}3(Y^aY^b\nd_bU_a)\si-\frac{3}{2}\rho Y^aY^b\nd_a\nd_b\rho+\frac{3}{2}(\tau+3\mu\rho)^2\\+9\rho^2\P_{ab}Y^aY^b+3\phi\si\rho-\frac{1}{2}(3\phi\rho-\si)^2\end{pmatrix}F^4\\
&+\begin{pmatrix}
\frac{1}{2}Y^aY^b\nd_a\nd_b\rho\si-\frac{185}{8}\rho^4-3(Y^aY^b\nd_bL_a)\rho-6\rho \si \P_{ab}Y^aY^b\\
+\phi \si(3\phi \rho-\si)+(3\mu \rho+\tau)(\rho \ell+\phi \tau)-(3\mu\rho +\tau)(Y^a\nd_a\si)
\end{pmatrix}F^2\\
&+\left(\frac{11}{2}\rho \si(\tau+3\rho \mu)-\frac{27}{2}\rho^2(\rho \ell+\phi \tau)-\frac{5}{2}\rho^2(Y^a\nd_a\si)\right)F\\
&+(Y^aY^b\nd_bL_a)\si-\frac{5}{2}\si \rho^3+\P_{ab}Y^aY^b \si^2-\frac{1}{2}(\rho\ell+\phi \tau)^2\\
&-\frac{1}{2}\phi^2\si^2-(\rho \ell+\tau \phi)(Y^a\nd_a\si)\\
=&0,
\end{align*}
and we clear denominators to obtain $P_2(F):=P_0(F)\tilde P_2(F)=0$. 
Finally the third constraint is given by
\begin{align*}
P_3(F)=&-6\tau F^6+18 \rho^2 F^5+(3Y^b\nd_b\si+24(\rho \ell+\tau \phi)-6\si \mu)F^4+13\si \rho F^3\\
&+\begin{pmatrix}(3\phi-\frac{\si}{\rho})Y^b\nd_b\si+30 \mu \phi \si+30 \phi \rho \ell+ 30 \phi^2 \tau\\-(3\mu +\frac{\tau}{\rho})(U^b\nd_b\si)-10\si \ell+3 \rho \epsilon^{ab}\nd_bL_a\end{pmatrix}F^2\\
&+\left(25\phi \si\rho-\frac{5\rho}{2}(U^b\nd_b\si)-8\si^2\right)F\\
&-\frac{\phi \si}{\rho}Y^b\nd_b\si-(U^b\nd_b\si)(\ell+\frac{\phi \tau}{\rho})-(\epsilon^{ab}\nd_bL_a)\si\\
=&0.
\end{align*} 
We have derived the polynomial constraints $P_1(F)=P_2(F)=P_3(F)=0$ for (\ref{sfMEW}) to hold under the assumption that the generic condition $P_0(F)\neq 0$ holds on $M^2$.
The polynomials $P_1(t)$, $P_2(t)$ and $P_3(t)$ in Theorem \ref{polycon} are obtained by replacing $F$ with the indeterminate $t$ (the polynomial $P_2(t)$ is given by $P_0(t)\tilde P_2(t)$).

\section{The case where $P_0(F)=0$}\label{mabproof}
In this section, we examine the case where $P_0(F)=\si-3\rho F^2=0$ and (\ref{sfMEW}) both hold on the M\"obius surface $(M^2,[g],[\P])$ and show that they imply the vanishing of a symmetric tensor $M_{ab}$ constructed out of invariants of the M\"obius structure. Under the assumption $P_0(F)=0$, we obtain from (\ref{abs}) that
\[
0=\rho \ell+\frac{5}{2}\rho^2 F+3\rho \mu F^2+\tau F^2+\tau \phi. 
\]
Substituting $F^2=\frac{\si}{3\rho}$ gives
\begin{align*}
0=\rho \ell+\frac{5}{2}\rho^2 F+\mu \si+\frac{\tau \si}{3 \rho}+\tau \phi,
\end{align*}
which upon rearranging gives
\begin{equation}\label{forf}
F=-\frac{2}{5}\left(\frac{\rho \ell+\mu \si+\frac{\tau \si}{3\rho}+\tau \phi}{\rho^2}\right).
\end{equation}
Also, the first constraint (\ref{mewcon2}) becomes
\begin{align*}
\al_aU^a=-(F^2+\phi)=-\left(\frac{\si}{3\rho}+\phi\right)=-m,
\end{align*}
where 
\[
m:=\frac{\si}{3\rho}+\phi.
\]
Differentiating once more and using (\ref{mewf1}) we obtain
\begin{align*}
\al_aW^a=&\al_a(Y^c\nd_cU^a+\phi Y^a-3\mu U^a)\\
=&-Y^c\nd_c m +\frac{1}{2}\rho F+\P_{ca}U^cY^a+(\phi-m)(\al_cY^c)+3\mu m. 
\end{align*}
Let 
\[
\psi=3\mu m+\P_{ca}U^cY^a-Y^c\nd_c m. 
\]
We obtain the following expression for $\al_a$:
\begin{align}\label{alphaf}
\al_a=\frac{Y_a}{\si}\left(\frac{1}{2}\rho F+\psi+(\phi-m)\al_cY^c\right)-\frac{\epsilon_{ab}W^b}{\si}m.
\end{align}
Contracting with $Y^a$ on both sides, we obtain
 \begin{align*}
\al_aY^a=&\frac{\rho}{\si}\left(\frac{1}{2}\rho F+\psi+(\phi-m)\al_cY^c\right)+\frac{\tau m}{\si}\\
=&\frac{\rho^2}{2\si} F+\frac{\psi \rho}{\si}-\frac{1}{3}\al_cY^c+\frac{\tau m}{\si}, 
\end{align*}
from which we obtain
\begin{align}\label{alphaff}
\al_aY^a=&\frac{3\rho^2}{8\si} F+\frac{3(\psi \rho+\tau m)}{4\si}. 
\end{align}
Substituting (\ref{alphaff}) into (\ref{alphaf}) now gives
\begin{align*}
\al_a=\left(\frac{3 \rho}{8 \si}F+\frac{3(\psi \rho+\tau m)}{4\si \rho}\right)Y_a-\frac{m U_a}{\rho}
\end{align*}
and a further substitution of (\ref{forf}) gives
\begin{align*}
\al_a=&\left(-\frac{3}{20}\left(\frac{\ell}{\si}+
\frac{\mu}{\rho}+
\frac{\tau}{3\rho^2}+
\frac{\tau \phi}{\rho \si}\right)+\frac{3(\psi \rho+\tau m)}{4\si \rho}\right)Y_a-\frac{m U_a}{\rho}.
\end{align*}
Defining the quantity $k$ by
\[
k:=-\frac{3\rho }{20}\left(\frac{\ell}{\si}+
\frac{\mu}{\rho}+
\frac{\tau}{3\rho^2}+
\frac{\tau \phi}{\rho \si}\right)+\frac{3(\psi \rho+\tau m)}{4\si}
\]
then gives 
\begin{align}\label{malpha}
\al_a=\frac{k Y_a}{\rho}-\frac{m U_a}{\rho}.
\end{align}
Let
\begin{equation}\label{defm}
M_{ab}:=\nd_{(a}\al_{b)}+\al_a\al_b+\P_{ab}-\frac{\al_c\al^c}{2}g_{ab}
\end{equation}
for $\al_a$ given by (\ref{malpha}). Hence for M\"obius structures with $P_0(F)=0$ and (\ref{sfMEW}) both holding, we necessarily must have $\al_a$ given by (\ref{malpha}), and the tensor $M_{ab}$ given by (\ref{defm}) automatically vanishes. Conversely for M\"obius structures with $M_{ab}=0$, taking $\al_a$ to be as given in (\ref{malpha}), we obtain a solution to (\ref{sfMEW}). This is the result of Theorem \ref{mab}. 

\section{Main Theorems}
The symmetric tensor $M_{ab}$ is an invariant of the M\"obius structure that can be used to distinguish M\"obius surfaces where sf-MEW holds with $P_0(F)=0$ from those with $P_0(F) \neq 0$. 
For the formulation of Theorem \ref{polycon}, we therefore assume that $M_{ab} \neq 0$, which allows us to work locally in an open set $U \subset M^2$ where $P_0(F)=\si-3\rho F^2 \neq 0$.
\begin{theorem}\label{polycon}
Let $(M^2,[g],[\P])$ be a M\"obius surface, with $M_{ab} \neq 0$. Suppose the surface locally admits a solution to (\ref{sfmew0}). Then there exist polynomials $P_1(t)$, $P_2(t)$, $P_3(t)$ in a single variable $t$ with coefficients given by conformal invariants of the M\"obius structure such that when $t=F$, where $F=\epsilon^{ab}F_{ab}$,
\[
P_1(F)=P_2(F)=P_3(F)=0
\]  
must hold.
\end{theorem}
We have explicitly computed the polynomial constraints $P_1(F)=P_2(F)=P_3(F)=0$ in section \ref{polyconproof}. 
For any polynomials $P(t)$, $Q(t)$ in a single variable $t$, let $\mbox{Res}(P(t),Q(t))$ denote the resultant of $P(t)$ and $Q(t)$. $\mbox{Res}(P(t),Q(t))=0$ is necessary and sufficient for $P(t)$ and $Q(t)$ to share a common root. As a corollary, we obtain local obstructions for there to be solutions of (\ref{sfmew0}).
\begin{corollary}\label{polyobs}
Suppose the M\"obius surface $(M^2,[g],[\P])$ has $M_{ab} \neq 0$ and locally admits a solution to (\ref{sfmew0}). Then the following conformal invariants of the M\"obius structure have to vanish: 
\[
\mbox{Res}(P_1(t),P_2(t))=\mbox{Res}(P_1(t),P_3(t))=\mbox{Res}(P_2(t),P_3(t))=0.
\]  
\end{corollary}
For the case when $M_{ab}=0$, we have
\begin{theorem}\label{mab}
Let $(M^2,[g],[\P])$ be a M\"obius surface with $M_{ab}=0$. Then the surface locally admits a solution to (\ref{sfmew0}), with $\al_a$ given by (\ref{malpha}).
\end{theorem}
However, the author does not know of any examples of non-flat M\"obius surfaces for which $M_{ab}=0$.

\section{Examples}\label{examples}
In this section we give examples of three different M\"obius structures on the Euclidean plane $\R^2$. The first example has non-vanishing obstruction, while the second and third have vanishing obstructions. In the third example we show that the M\"obius structure does not admit a real solution to (\ref{sfmew0}) despite having vanishing obstructions.  
Since $\rho=Y_aY^a>0$ for non-flat M\"obius surfaces, computing the discriminant of $P_0(F)$ gives $12\rho \si > 0$, or $\si >0$ for there to be solutions of $P_0(F)=0$. (For $\si=0$, it would mean $-3\rho F^2=0$ which cannot happen if $\rho>0$ and $F \neq 0$.) 
In the first $2$ examples we have $\si \leq 0$, so that $P_0(F)=\si-3\rho F^2 \neq 0$ and we do not need to compute the M\"obius invariant $M_{ab}$. The last example has $\si>0$, but a simple argument eliminates the possibility that $P_0(F)=0$ can hold.  
\subsection{Example with non-vanishing obstruction}
The first example will be the M\"obius structure given by $\P_{ab}=x_{(a}\epsilon_{b)c}x^c$ on $\R^2$ with the flat metric $\de_{ab}$. On $\R^2$ we have $K=0$. Here $x^a$ are standard local coordinates in $\R^2$ so that $\partial_ax_b=\de_{ab}$. Let $r=x_ax^a$.
For this M\"obius structure, we have
\begin{align*}
P_1(t)=&256r^2\left(\frac{63}{2}t^8+56t^4-640rt^3+672r^2t^2+320r^3t+32r^4\right),\\
P_2(t)=&256r^2\left(-\frac{9}{2}t^8+56t^4-7392r^2t^2+1472r^3t+288r^4\right),\\
P_3(t)=&-768r t^4( t^2-6rt-4r^2).
\end{align*}
Using MAPLE, we find 
\begin{align*}
\mbox{Res}(P_1(t),P_3(t))=2^{142}\cdot3^{10}r^{44}(2^4\cdot3^{2}\cdot7^2r^{8}+ 2^2 \cdot 7^2 \cdot 59 \cdot 251r^{4}-3^2 \cdot 131)
\end{align*}
which is non-zero for general $r$. Similarly, the obstructions given by $\mbox{Res}(P_2(t),P_3(t))$ and $\mbox{Res}(P_1(t),P_2(t))$ do not vanish on any open set. We conclude that the M\"obius structure for this example admits no local solution to (\ref{sfMEW}).    

\subsection{Example with vanishing obstruction (and solution to (\ref{sfMEW}))}
Consider the M\"obius structure given by $\P_{ab}=x_ax_b-\frac{1}{2}\de_{ab}x_cx^c$ on $\R^2$ with the flat metric $\de_{ab}$. Again we have $K=0$. It can be verified that $\al_a=\pm\epsilon_{ab}x^b$ is a solution to (\ref{sfMEW}), since 
\[
\partial_a\al_b=\partial_{a}(\pm\epsilon_{bc}x^c)=\pm\epsilon_{ba}=\mp\epsilon_{ab}, 
\] 
and therefore, using $\epsilon_{ac}\epsilon_{bd}=\de_{ab}\de_{cd}-\de_{ad}\de_{cb}$, we find
\begin{align*}
\partial_{(a}\al_{b)}+\al_a\al_b+\P_{ab}-\frac{\al_c\al^c}{2}\de_{ab}
=&0+\epsilon_{ac}x^c\epsilon_{bd}x^d+(x_ax_b-\frac{1}{2}\de_{ab}x_cx^c)-\frac{x_cx^c}{2}\de_{ab}\\
=&\de_{ab}x_cx^c-x_ax_b+x_ax_b-\frac{1}{2}\de_{ab}x_cx^c-\frac{x_cx^c}{2}\de_{ab}\\
=&0.
\end{align*}
We have 
\[
F_{ab}=\partial_{[a}\al_{b]}=\pm\epsilon_{[ba]}=\mp\epsilon_{ab}, 
\]
so that $F=\mp2$ depending on the sign of $\al_a$ chosen. Computing $P_1(t)$ gives
\begin{align*}
P_1(t)=&\rho^2(t^2-4)\left(\frac{63}{2}t^6+222t^4+(512-\frac{9\rho^2}{32})t^2+(384+\frac{\rho^2}{2})\right).
\end{align*}
Computing $\tilde P_2(t)$ gives
\begin{align*}
\tilde P_2(t)
=&\left(\frac{8+15 t^2}{8+3 t^2}\right)\left(\frac{25}{4}\rho^4 t^2\right)-\frac{9}{2}\rho^2t^8+48\rho^2t^6+(136\rho^2-\frac{9\rho^4}{32})t^4\\
&+(-640\rho^2-\frac{197}{8}\rho^4)t^2-1536\rho^2+18\rho^4.
\end{align*}
Multiplying throughout by $P_0(t)=-\rho(8+3t^2)$ and expanding the terms on the right, we obtain
\begin{align*}
P_2(t)=&-\rho(8+3t^2)\tilde P_2(t)\\
=&\frac{27}{2}\rho^3(t^2-4)\left(t^8-4t^6+(\frac{\rho^2}{16}-\frac{224}{3})t^4-(\frac{6400}{27}+\frac{19}{18}\rho^2)t^2+(\frac{8}{3}\rho^2-\frac{2048}{9})\right).
\end{align*}
The third polynomial $P_3(t)$ for this example is
\begin{align*}
P_3(t)=&2\rho^2t(9t^2-16)(t^2-4).
\end{align*}
It therefore can be seen that $P_1(t)$, $P_2(t)$ and $P_3(t)$ share a common root $t^2-4$, attained when $t=F$. The local obstructions vanish for this example. 

\subsection{Example with vanishing obstructions but does not admit a real solution to (\ref{sfMEW})}
In this example we show that obstructions can vanish for a particular M\"obius structure yet it does not admit a real solution to (\ref{sfMEW}).
Consider the M\"obius structure given by $\P_{ab}=\frac{1}{2}\de_{ab}x_cx^c-x_ax_b$ on $\R^2$ with the flat metric $\de_{ab}$. For this M\"obius structure, we have $\si=8\rho$, and $\mu=0$, $\ell=0$, $\tau=0$. The polynomial constraint $P_0(F)=0$ cannot hold because equation (\ref{forf}) yields $F=0$, which cannot happen if the M\"obius structure is not flat. Proceeding to compute $P_1(t)$ for this example, we find that
\begin{align*}
P_1(t)
=&\rho^2(t^2+4)\left(\frac{63}{2}t^6-222t^4+(512+\frac{9\rho^2}{32})t^2+(\frac{\rho^2}{2}-384)\right),
\end{align*}
and we also have
\begin{align*}
P_2(t)=&(8-3t^2)\tilde P_2(t)\\
=&\frac{27}{2}\rho^2(t^2+4)\left(t^8+4t^6-(\frac{\rho^2}{16}+\frac{224}{3})t^4+(\frac{6400}{27}-\frac{19}{18}\rho^2)t^2-(\frac{8}{3}\rho^2+\frac{2048}{9})\right)
\end{align*}
and 
\begin{align*}
P_3(t)=&2\rho^2t(9t^2+16)(t^2+4).
\end{align*}
The three polynomials $P_1(t)$, $P_2(t)$ and $P_3(t)$ all share a common factor $t^2+4$. Hence the local obstructions all vanish for this example. We now deduce that the M\"obius structure does not admit any real solution to (\ref{sfMEW}) by showing that these three polynomials do not share any other common factors besides $t^2+4$. 
Indeed, dividing $P_1(t)$, $P_2(t)$ and $P_3(t)$ by their common factors $\rho^2(t^2+4)$, we obtain
\begin{align*}
S_1(t)=&\frac{63}{2}t^6-222t^4+(512+\frac{9\rho^2}{32})t^2+(\frac{\rho^2}{2}-384),\\
S_2(t)=&\frac{27}{2}\left(t^8+4t^6-(\frac{\rho^2}{16}+\frac{224}{3})t^4+(\frac{6400}{27}-\frac{19}{18}\rho^2)t^2-(\frac{8}{3}\rho^2+\frac{2048}{9})\right),\\
S_3(t)=&2t(9t^2+16).
\end{align*}
We find that
\begin{align*}
\mbox{Res}(S_1(t)&,S_2(t))\\
=&\frac{1076168025}{67108864}\rho^8(243\rho^6+12704256\rho^4+131135897600\rho^2+251658240000)^2
\end{align*}
and
\begin{align*}
\mbox{Res}(S_1(t),S_3(t))=&80289792000000\rho^2-61662560256000000,\\
\mbox{Res}(S_2(t),S_3(t))=&-3583180800(73600+81\rho^2)^2(3\rho^2+256),
\end{align*}
and so the three polynomials share no other common factors. 
However, the M\"obius structure does admit a complex solution to (\ref{sfMEW}). It can be verified that $\al_a=\pm i \epsilon_{ab}x^b$ is a solution to (\ref{sfMEW}), 
since 
\[
\partial_a\al_b=\partial_{a}(\pm i\epsilon_{bc}x^c)=\pm i\epsilon_{ba}=\mp i\epsilon_{ab}, 
\] 
and therefore
\begin{align*}
\partial_{(a}\al_{b)}+\al_a\al_b+\P_{ab}-\frac{\al_c\al^c}{2}\de_{ab}
=&0-\epsilon_{ac}x^c\epsilon_{bd}x^d+(\frac{1}{2}\de_{ab}x_cx^c-x_ax_b)+\frac{x_cx^c}{2}\de_{ab}\\
=&-\de_{ab}x_cx^c+x_ax_b-x_ax_b+\frac{1}{2}\de_{ab}x_cx^c+\frac{x_cx^c}{2}\de_{ab}\\
=&0.
\end{align*}
We have 
\[
F_{ab}=\partial_{[a}\al_{b]}=\pm i\epsilon_{[ba]}=\mp i\epsilon_{ab}, 
\]
so that $F=\mp2i$ depending on the sign of $\al_a$ chosen.

\end{document}